\documentclass[11pt,reqno]{amsart}

\usepackage{amssymb,amsmath,epsfig}
\usepackage{amsfonts}

\newcommand{\Z}{\mathbb{Z}}

\newcommand\smod[1]{(\operatorname{mod} #1)}
\renewcommand{\theta}{\vartheta}
\newtheorem{theorem}{Theorem}[section]

\newtheorem{lemma}[theorem]{Lemma}

\newtheorem{definition}[theorem]{Definition}
\newtheorem*{theorem*}{Theorem}

\numberwithin{equation}{section}

\title{The tenth order mock theta functions revisited}

\date{\today}
\author{Sander Zwegers}
\address{School of Mathematical Sciences, University College Dublin, Belfield, Dublin 4, Ireland}
\email{sander.zwegers@ucd.ie}


\begin{document}
\begin{abstract}
In this paper we consider the first four of the eight identities between the tenth order mock theta functions, found in Ramanujan's lost notebook. These were originally proved by Choi. Here we give an alternative proof.
\end{abstract}
\maketitle
\section{Introduction and statement of results}
The mock theta functions were the subject of the last letter that Ramanujan wrote to Hardy, shortly before his death. In this letter he gives a list of mock theta functions of ``order'' three, five and seven, together with identities that they satisfy.

In Ramanujan's lost notebook \cite{Ram} we can find identities between further examples of mock theta functions. For example, eight identities between ``tenth order'' mock theta functions are given. These identities were established by Choi in a series of papers, see \cite{Choi1,Choi2,Choi3,Choi4}. Each of these papers is devoted to two of Ramanujan's identities. The aim of this paper is to give shorter proofs for the first four identities (proved by Choi in \cite{Choi1,Choi2}). The identities involve the four mock theta functions $\phi$, $\psi$, $X$ and $\chi$. Ramanujan defined these in terms of particular $q$-hypergeometric series. For example
\begin{equation*}
\phi (q) := \frac{1}{1-q} + \frac{q}{(1-q)(1-q^3)} + \frac{q^3}{(1-q)(1-q^3)(1-q^5)} + \frac{q^6}{(1-q)(1-q^3)(1-q^5)(1-q^7)} +\ldots
\end{equation*}
See \cite{Choi1} for the definitions of the other three. For the first part of Choi's proof, he finds identities for the four functions, which relate them to indefinite theta functions. We will cheat a bit and use these identities as our starting point. In fact, for the purpose of this paper, we will use these as the definitions of the four functions. In section 2.2 (pp.\ 507--513) of \cite{Choi1} he finds the following, slightly rewritten identities
\begin{equation}\label{indef1}
\begin{split}
\phi(q) &= \frac{1}{\sum_{n\in\Z} (-1)^n q^{n^2} } \left( \sum_{r,s \geq 0} - \sum_{r,s <0} \right) (-1)^{r+s} q^{r^2+3rs+s^2 +r+s},\\
\psi(q) &= \frac{1}{\sum_{n\in\Z} (-1)^n q^{n^2} } \left( \sum_{r,s \geq 0} - \sum_{r,s <0} \right) (-1)^{r+s+1} q^{r^2+3rs+s^2 +3r+3s+2},
\end{split}
\end{equation}
and in section 2.2 (pp.\ 190--197) of \cite{Choi2} he finds
\begin{equation}\label{indef2}
\begin{split}
X(q) &= \frac{2}{\sum_{n\in\Z} q^{\frac{1}{2}n(n+1)} } \left( \sum_{r,s \geq 0} - \sum_{r,s <0} \right) q^{2r^2+6rs+2s^2 +r+s},\\
\chi(q) &= 2- \frac{2q}{\sum_{n\in\Z} q^{\frac{1}{2}n(n+1)} } \left( \sum_{r,s \geq 0} - \sum_{r,s <0} \right) q^{2r^2+6rs+2s^2 -3r-3s}.
\end{split}
\end{equation}
For the rewriting we have used the classical identity $\sum_{n\in\Z} q^{\frac{1}{2}n(n+1)} = 2 (q^2;q^2)_\infty^2/(q;q)_\infty$, where we use the standard notation $(x;q)_\infty = \prod_{n=0}^\infty (1-xq^n)$.

The first four identities given by Ramanujan then are (see \cite[pp.\ 181--182]{Choi2}, again slightly rewritten)
\begin{equation}\label{id}
\begin{split}
q^{2/3} \phi (q^3) - \frac{\psi(\omega q^{1/3}) - \psi(\omega^2 q^{1/3})}{\omega-\omega^2} &= - q^{1/3} \frac{\sum_{n\in\Z} (-1)^n q^{n^2/3}}{\sum_{n\in\Z} (-1)^n q^{n^2}}\frac{\sum_{n\in\Z} (-1)^n q^{5n^2/2 +3n/2}}{(q;q^2)_\infty},\\
q^{-2/3} \psi (q^3) + \frac{\omega\phi(\omega q^{1/3}) - \omega^2 \phi(\omega^2 q^{1/3})}{\omega-\omega^2} &= \frac{\sum_{n\in\Z} (-1)^n q^{n^2/3}}{\sum_{n\in\Z} (-1)^n q^{n^2}}\frac{\sum_{n\in\Z} (-1)^n q^{5n^2/2 +n/2}}{(q;q^2)_\infty},\\
X(q^3) - \frac{\omega\chi(\omega q^{1/3}) - \omega^2 \chi(\omega^2 q^{1/3})}{\omega-\omega^2} &= \frac{\sum_{n\in\Z} q^{n(n+1)/6}}{\sum_{n\in\Z} q^{n(n+1)/2}}\frac{\sum_{n\in\Z} (-1)^n q^{5n^2 +n}}{(-q;q)_\infty},\\
\chi(q^3) + q^{2/3} \frac{X(\omega q^{1/3}) - X(\omega^2 q^{1/3})}{\omega-\omega^2} &= -q\frac{\sum_{n\in\Z} q^{n(n+1)/6}}{\sum_{n\in\Z} q^{n(n+1)/2}}\frac{\sum_{n\in\Z} (-1)^n q^{5n^2 +3n}}{(-q;q)_\infty},
\end{split}
\end{equation}
where $\omega$ is a primitive cube root of unity. Note that the exponent 7 in the second formula of (1.5) in \cite{Choi2} should be a 9.

Recent developments have shed new light on Ramanujan's mock theta functions. It turns out that they can be seen as the holomorphic parts of certain weight 1/2 harmonic Maass forms (for example see \cite{Ono,Zag}). One application of this theory is that it now becomes as easy to prove identities between mock theta functions as it is for modular forms. The proof is essentially automatic: using results from \cite{zwe2} (see also \cite{zwe1}) we can find a non-holomorphic correction term to make the mock theta functions transform like modular forms and we can find the explicit modular transformation properties; we then check that the correction terms cancel and so the proof is reduced to proving an identity between (weakly) holomorphic modular forms, which can be done with classical means. This paper started as an attempt to do this explicitly for the identities between the tenth order mock theta functions. Although the method is straightforward, to work out all the details is a bit tedious. To simplify things a bit, we first write the identities, using elementary methods, in a more suitable form. For this we introduce
\begin{definition}
For $r,s\in\Z$ we define
\begin{equation*}
\begin{split}
\rho_{r,s} &:= \begin{cases}
1& \text{if}\ r,s\geq 0,\\
-1 & \text{if}\ r,s<0,\\
0 & \text{otherwise}.
\end{cases}\\
\delta(r) &:= \begin{cases}
1 & \text{if}\ r\equiv 0 \smod{3},\\
0 & \text{otherwise}.
\end{cases}
\end{split}
\end{equation*}
\end{definition}
We will see that the identities in \eqref{id} can then be rewritten as
\begin{theorem}\label{theo}
We have
\begin{equation*}
\begin{split}
\sum_{k,l,r,s\in\Z} \rho_{r,s}& (-1)^{k+l+r+s} \bigl(\delta(k) - \delta (r) \bigr) \bigl(\delta(l) - \delta(s)\bigr)\ q^{(k^2+l^2+r^2+3rs+s^2 +3r+3s +1)/3} \\
&=-(q;q)_\infty \left(\sum_{n\in\Z} (-1)^n q^{n^2} \right)^2 \sum_{n\in\Z} (-1)^n q^{5n^2/2+ 3n/2},\\
\sum_{k,l,r,s\in\Z} \rho_{r,s}& (-1)^{k+l+r+s} \bigl(\delta(k) - \delta (r-1) \bigr) \bigl(\delta(l) - \delta(s-1)\bigr)\ q^{(k^2+l^2+r^2+3rs+s^2 +r+s)/3} \\
&=(q;q)_\infty \left(\sum_{n\in\Z} (-1)^n q^{n^2} \right)^2 \sum_{n\in\Z} (-1)^n q^{5n^2/2+ n/2},\\
\sum_{k,l,r,s\in\Z} \rho_{r,s}& \bigl(\delta(k-1) - \delta (r) \bigr) \bigl(\delta(l-1) - \delta(s)\bigr)\  q^{(k(k+1)/2 + l(l+1)/2+2r^2+6rs +2s^2 +3r+3s)/3}\\
&=(q^2;q^2)_\infty \left(\sum_{n\in\Z} q^{\frac{1}{2}n(n+1)} \right)^2 \sum_{n\in\Z} (-1)^n q^{5n^2+ n},\\
\sum_{k,l,r,s\in\Z} \rho_{r,s}& \bigl(\delta(k-1) - \delta (r+1) \bigr) \bigl(\delta(l-1) - \delta(s+1)\bigr) q^{(k(k+1)/2 + l(l+1)/2+2r^2+6rs +2s^2 +r+s-2)/3}\\
&=(q^2;q^2)_\infty \left(\sum_{n\in\Z} q^{\frac{1}{2}n(n+1)} \right)^2 \sum_{n\in\Z} (-1)^n q^{5n^2+ 3n}.
\end{split}
\end{equation*}
\end{theorem}
Using results from the theory of indefinite quadratic forms of type $(r,1)$, see \cite{zwe2} and also \cite{zwe1}, we can show that the left hand sides are all holomorphic modular forms and we can explicitly get the modular transformation properties. The advantage of this rewriting is that it simplifies the proof, because we don't have to explicitly compute the non-holomorphic correction terms in the original equations. Although the methods to prove Theorem \ref{theo} are again straightforward, to work out all the details is still a bit tedious. Therefore, we will use the so-called constant term method to prove Theorem \ref{theo} directly, instead.

The outline of this paper is as follows: in the next section we will show that the equations \eqref{id} are equivalent to the ones in Theorem \ref{theo}. In Section 3 we will then prove Theorem \ref{theo}.

\section{Rewriting the equations}
Note that throughout the paper, if not mentioned explicitly, the summation indices run through all integers, or through all integers satisfying the conditions listed under the summation sign.
 
We start with a lemma.
\begin{lemma}
Define
\begin{equation*}
\theta_0 (q) := \sum_{n\in\Z} (-1)^n q^{n^2} \qquad \text{and}\qquad 
\theta_1 (q) := \frac{1}{2} \sum_{n\in\Z} q^{n(n+1)/2},
\end{equation*}
then
\begin{equation*}
\begin{split}
\theta_j (q^{1/3}) \theta_j (\omega q^{1/3}) \theta_j (\omega^2 q^{1/3}) \theta_j (q^3)&= \theta_j(q)^4,
\end{split}
\end{equation*}
for $j=0,1$.
\end{lemma}
\begin{proof}
Checking term by term and using $(1-x)(1-\omega x)(1-\omega^2 x)=(1-x^3)$, we can easily verify that
\begin{equation*}
(q^{1/3};q^{1/3})_\infty (\omega q^{1/3}; \omega q^{1/3})_\infty (\omega^2 q^{1/3}; \omega^2 q^{1/3})_\infty (q^3;q^3)_\infty = (q;q)_\infty^4.
\end{equation*}
The result then follows from the classical identities
\begin{equation}\label{theta}
\begin{split}
\theta_0 (q) &= \frac{(q;q)_\infty^2}{ (q^2;q^2)_\infty},\\
\theta_1 (q) &= \frac{(q^2;q^2)_\infty^2}{(q;q)_\infty}.
\end{split}
\end{equation}
\end{proof}
We start with the first equation in \eqref{id}. If we multiply both sides by 
\begin{equation*}
q^{-1/3} \theta_0(\omega q^{1/3}) \theta_0(\omega^2 q^{1/3}) \theta_0 (q^3),
\end{equation*}
then using the lemma with $j=0$, equation \eqref{theta} and $(q;q^2)_\infty = (q;q)_\infty / (q^2;q^2)_\infty$, we see that the right hand side becomes the right hand side of the first equation in Theorem \ref{theo}. Using \eqref{indef1} we see that the left hand side becomes
\begin{equation}\label{re1}
\begin{split}
q^{1/3} \theta_0 (\omega q^{1/3}) \theta_0 (\omega^2 q^{1/3}) &\sum_{r,s} \rho_{r,s} (-1)^{r+s} q^{3r^2+9rs +3s^2 +3r+3s}\\
- q^{-1/3} \frac{\theta_0 (q^3)}{\omega - \omega^2} &\left\{\theta_0( \omega^2 q^{1/3}) \sum_{r,s} \rho_{r,s} (-1)^{r+s+1} \left(\omega q^{1/3}\right)^{r^2+3rs+s^2+3r+3s+2}\right. \\
&- \left.\theta_0( \omega q^{1/3}) \sum_{r,s} \rho_{r,s} (-1)^{r+s+1} \left(\omega^2 q^{1/3}\right)^{r^2+3rs+s^2+3r+3s+2}\right\}.
\end{split}
\end{equation}
Using $\theta_0 (q^{1/3}) = \theta_{0,0}(q) + \theta_{0,1} (q)+ \theta_{0,2}(q)$ with $\theta_{0,j} = \sum_{n\equiv j\pmod{3}} (-1)^n q^{n^2/3}$ we find
\begin{equation*}
\begin{split}
\theta_0 (\omega q^{1/3}) \theta_0 (\omega^2 q^{1/3})&= \bigl(\theta_{0,0}(q) + \omega\theta_{0,1} (q)+ \omega\theta_{0,2}(q)\bigr)\bigl(\theta_{0,0}(q) + \omega^2\theta_{0,1} (q)+ \omega^2\theta_{0,2}(q)\bigr)\\
&= \theta_{0,0}(q)^2 -\theta_{0,0}(q)\theta_{0,1}(q)-\theta_{0,0}(q)\theta_{0,2}(q)+\theta_{0,1}(q)^2+2\theta_{0,1}(q)\theta_{0,2}(q)+\theta_{0,2}(q)^2\\
&= \sum_{k,l\in\Z} p(k,l) (-1)^{k+l} q^{(k^2+l^2)/3},
\end{split}
\end{equation*}
where 
\begin{equation*}
p(k,l)= \begin{cases}
1& \text{if}\ (k,l)\equiv (0,0), (1,1), (1,2), (2,1), (2,2) \pmod{3},\\
-1 & \text{if}\ (k,l)\equiv (0,1), (0,2) \pmod{3},\\
0 & \text{if}\ (k,l)\equiv (1,0), (2,0) \pmod{3}.
\end{cases}
\end{equation*}
Hence the first term in \eqref{re1} can be rewritten as
\begin{equation}\label{re2}
\sum_{k,l,r,s} \rho_{r,s}\ \delta(r) \delta(s) p(k,l) (-1)^{k+l+r+s} q^{(k^2+l^2+r^2+3rs +s^2 +3r+3s+1)/3}.
\end{equation}
The second term in \eqref{re1} we rewrite using $\frac{\omega^k - \omega^{2k}}{\omega - \omega^2}= \chi_3 (k)$, with 
\begin{equation*}
\chi_3(k) =\begin{cases}
0 & \text{if}\ k\equiv 0 \pmod{3},\\
1 & \text{if}\ k\equiv 1 \pmod{3},\\
-1 & \text{if}\ k\equiv 2 \pmod{3},\\
\end{cases}
\end{equation*}
as
\begin{equation}\label{re3}
\begin{split}
\theta_0 (q^3) &\sum_{l,r,s} \rho_{r,s}\ \chi_3(2l^2+ r^2+3rs+s^2+3r+3s+2) (-1)^{l+r+s} q^{(l^2+r^2+3rs +s^2 +3r+3s+1)/3}\\
&=\sum_{k,l,r,s} \rho_{r,s}\ \delta(k) \chi_3(r^2+s^2-l^2-1) (-1)^{k+l+r+s} q^{(k^2+l^2+r^2+3rs +s^2 +3r+3s+1)/3}.
\end{split}
\end{equation}
If we combine the equations \eqref{re2} and \eqref{re3}, we see that the expression in \eqref{re1} equals
\begin{equation*}
\sum_{k,l,r,s} \rho_{r,s}\ \bigl(\delta(r) \delta(s) p(k,l)+\delta(k) \chi_3(r^2+s^2-l^2-1)\bigr) (-1)^{k+l+r+s} q^{(k^2+l^2+r^2+3rs +s^2 +3r+3s+1)/3}.
\end{equation*}
If we use that
\begin{equation}\label{del}
\begin{split}
\delta(r) \delta(s) p(k,l)+&\delta(k) \chi_3(r^2+s^2-l^2-1) \\
&= \bigl(\delta(k) - \delta (r) \bigr) \bigl(\delta(l) - \delta(s)\bigr) + \delta(r) \bigl(\delta(s)-1\bigr)\bigl(\delta(k)-\delta(l)\bigr),
\end{split}
\end{equation}
which we can easily verify, since there are only 81 cases to check, we get the left hand side of the first equation in Theorem \ref{theo}. Note that the term $\delta(r) (\delta(s)-1)(\delta(k)-\delta(l)) $ cancels in the sum, because it is anti-symmetric in $k$ and $l$.

Hence we have established that the first equation in \eqref{id} is equivalent to the first equation in Theorem \ref{theo}. The arguments for the other equations are very similar, so we omit some of the details.

For the second equation in \eqref{id} we multiply both sides by 
\begin{equation*}
\theta_0(\omega q^{1/3}) \theta_0(\omega^2 q^{1/3}) \theta_0 (q^3).
\end{equation*}
The right hand side is easily seen to become the right hand side of the second equation in Theorem \ref{theo}. Using \eqref{indef1} we can again rewrite the left hand side. We get the sum of two terms. The first is
\begin{equation*}
\begin{split}
q^{-2/3} \theta_0(\omega q^{1/3}) \theta_0(\omega^2 q^{1/3}) &\theta_0 (q^3) \psi (q^3) \\
&=  \sum_{k,l,r,s} \rho_{r,s}\ \delta(r) \delta(s) p(k,l) (-1)^{k+l+r+s+1} q^{(k^2+l^2+r^2+3rs +s^2 +9r+9s+16)/3}\\
&= \sum_{k,l,r,s} \rho_{r,s}\ \delta(r+2) \delta(s+2) p(k,l) (-1)^{k+l+r+s} q^{(k^2+l^2+r^2+3rs +s^2 +r+s)/3},
\end{split}
\end{equation*}
were in the last step we replace $(r,s)$ by $(-r-2,-s-2)$ and use $\rho_{-r-2,-s-2}=-\rho_{r,s} -d(r+1)-d(s+1)$, where $d(n)$ is 1 if $n=0$ and 0 otherwise, and that $d(n+1)\delta(n+2)=0$ for all $n\in\Z$.

The second term is
\begin{equation*}
\begin{split}
\theta_0 (q^3) &\sum_{l,r,s} \rho_{r,s}\ \chi_3(2l^2+ r^2+3rs+s^2+r+s+1) (-1)^{l+r+s} q^{(l^2+r^2+3rs +s^2 +r+s)/3}\\
&=\sum_{k,l,r,s} \rho_{r,s}\ \delta(k) \chi_3(r^2+s^2-l^2+r+s+1) (-1)^{k+l+r+s} q^{(k^2+l^2+r^2+3rs +s^2 +r+s)/3}.
\end{split}
\end{equation*}
If we combine these two terms and use that
\begin{equation*}
\begin{split}
\delta(r+2) \delta(s+2) p(k,l) + &\delta(k) \chi_3(r^2+s^2-l^2+r+s+1)\\
&=\bigl(\delta(k) - \delta (r-1) \bigr) \bigl(\delta(l) - \delta(s-1)\bigr) + \delta(r+2) \bigl(\delta(s+2)-1\bigr)\bigl(\delta(k)-\delta(l)\bigr),
\end{split}
\end{equation*}
which we get from \eqref{del} by replacing $(r,s)$ by $(r+2,s+2)$, we get the left hand side of the second equation in Theorem \ref{theo}.

For the third and fourth equation we first have to rewrite the identity for $\chi$ given in \eqref{indef2}. If we replace $(r,s)$ by $(-r,-s)$ and use that $\rho_{-r,-s} =- \rho_{r,s} +d(r)+d(s)$ we get
\begin{equation*}
\sum_{r,s} \rho_{r,s} q^{2r^2+6rs+2s^2 -3r-3s}= -\sum_{r,s} \rho_{r,s} q^{2r^2+6rs+2s^2 +3r+3s}+2\sum_n q^{2n^2+3n}.
\end{equation*}
Using
\begin{equation*}
2q \sum_n q^{2n^2+3n}= 2\sum_{m\ \text{odd}}q^{\frac{1}{2}m(m+1)} = \sum_m q^{\frac{1}{2}m(m+1)}
\end{equation*}
we then get from \eqref{indef2}
\begin{equation*}
\chi(q) = \frac{2}{\sum_{n\in\Z} q^{\frac{1}{2}n(n+1)} } \sum_{r,s} \rho_{r,s} q^{2r^2+6rs+2s^2 +3r+3s+1}.
\end{equation*}
We will use this identity instead of the one in \eqref{indef2}.

In the remaining two equations we multiply both sides by 
\begin{equation*}
4\theta_1(\omega q^{1/3}) \theta_1(\omega^2 q^{1/3}) \theta_1 (q^3) \qquad \text{and} \qquad -4q^{-1}\theta_1(\omega q^{1/3}) \theta_1(\omega^2 q^{1/3}) \theta_1 (q^3)
\end{equation*}
respectively. Using the lemma with $j=1$, equation \eqref{theta} and $(-q;q)_\infty = (q^2;q^2)_\infty / (q;q)_\infty$, we see that the right hand sides become the right hand sides of the third and fourth equation in Theorem \ref{theo}. To rewrite the left hand sides we use the equation $\theta_1 (q^{1/3}) = \theta_{1,0}(q) + \theta_{1,1} (q)+ \theta_{1,2}(q)$ with $\theta_{1,j} = \frac{1}{2}\sum_{n\equiv j\pmod{3}} q^{n(n+1)/6}$ to find
\begin{equation*}
\begin{split}
\theta_1 (\omega q^{1/3}) \theta_1 (\omega^2 q^{1/3})&= \bigl(\theta_{1,0}(q) + \omega\theta_{1,1} (q)+\theta_{1,2}(q)\bigr)\bigl(\theta_{1,0}(q) + \omega^2\theta_{1,1} (q)+ \theta_{1,2}(q)\bigr)\\
&= \theta_{1,0}(q)^2+2\theta_{1,0}(q)\theta_{1,2}(q)+\theta_{1,2}(q)^2 -\theta_{1,0}(q)\theta_{1,1}(q)-\theta_{1,1}(q)\theta_{1,2}(q)+\theta_{1,1}(q)^2\\
&= \frac{1}{4} \sum_{k,l\in\Z} p(k-1,l-1) q^{k(k+1)/6 + l(l+1)/6}.
\end{split}
\end{equation*}
The first term in the third equation then is
\begin{equation*}
\sum_{k,l,r,s} \rho_{r,s}\ \delta(r) \delta(s) p(k-1,l-1) q^{(k(k+1)/2+l(l+1)/2+2r^2+6rs +2s^2 +3r+3s)/3}
\end{equation*}
and the second
\begin{equation*}
\begin{split}
-2\theta_1 (q^3) &\sum_{l,r,s} \rho_{r,s}\ \chi_3(l^2+l+ 2r^2+6rs+2s^2+3r+3s+2) q^{(l(l+1)/2+2r^2+6rs +2s^2 +3r+3s+1)/3}\\
&=-\sum_{k,l,r,s} \rho_{r,s}\ \delta(k-1) \chi_3(-r^2-s^2+l^2+l-1) q^{(k(k+1)/2 + l(l+1)/2+2r^2+6rs +2s^2 +3r+3s)/3}.
\end{split}
\end{equation*}
If we now use
\begin{equation*}
\begin{split}
\delta(r) \delta(s) p(k-1,l-1)-&\delta(k-1) \chi_3(-r^2-s^2+l^2+l-1) \\
&= \bigl(\delta(k-1) - \delta (r) \bigr) \bigl(\delta(l-1) - \delta(s)\bigr) + \delta(r) \bigl(\delta(s)-1\bigr)\bigl(\delta(k-1)-\delta(l-1)\bigr),
\end{split}
\end{equation*}
which follows from \eqref{del} if we replace $(k,l)$ by $(k-1,l-1)$, we get the desired result.

After multiplication, the first term in the fourth equation is
\begin{equation*}
\begin{split}
-\sum_{k,l,r,s} &\rho_{r,s}\ \delta(r) \delta(s) p(k-1,l-1) q^{(k(k+1)/2+l(l+1)/2+2r^2+6rs +2s^2 +9r+9s+6)/3}\\
&= \sum_{k,l,r,s} \rho_{r,s}\ \delta(r+1) \delta(s+1) p(k-1,l-1) q^{(k(k+1)/2+l(l+1)/2+2r^2+6rs +2s^2 +r+s-2)/3},
\end{split}
\end{equation*}
were in the last step we replaced $(r,s)$ by $(-r-1,-s-1)$ and used that $\rho_{-r-1,-s-1}=-\rho_{r,s}$.
The second term is
\begin{equation*}
\begin{split}
-2q^{-1/3} &\theta_1 (q^3) \sum_{l,r,s} \rho_{r,s}\ \chi_3(l^2+l+ 2r^2+6rs+2s^2+r+s) q^{(l(l+1)/2+2r^2+6rs +2s^2 +r+s)/3}\\
&=-\sum_{k,l,r,s} \rho_{r,s}\ \delta(k-1) \chi_3(-r^2-s^2+l^2+l+r+s) q^{(k(k+1)/2 + l(l+1)/2+2r^2+6rs +2s^2 +r+s-2)/3}.
\end{split}
\end{equation*}
If we now use
\begin{equation*}
\begin{split}
\delta(r+1) &\delta(s+1) p(k-1,l-1)-\delta(k-1) \chi_3(-r^2-s^2+l^2+l+r+s) \\
& = \bigl(\delta(k-1) - \delta (r+1) \bigr) \bigl(\delta(l-1) - \delta(s+1)\bigr) - \delta(r+1) \bigl(\delta(s+1)-1\bigr)\bigl(\delta(k-1)-\delta(l-1)\bigr),
\end{split}
\end{equation*}
which follows from \eqref{del} by replacing $(k,l,r,s)$ by $(k-1,l-1,r+1,s+1)$, we get the desired result.
\section{Proof of Theorem \ref{theo}}
For the proof we need the following lemma.
\begin{lemma}\label{lem1}
Let 
\begin{equation}\label{jac}
\Theta(x;q) := \sum_{n\in\Z} (-1)^n q^{n(n-1)/2}x^n = (q;q)_\infty (x;q)_\infty (x^{-1} q;q)_\infty,
\end{equation}
where the last identity is the Jacobi triple product identity. Then
\begin{equation*}
\begin{split}
\sum_{k,l\in\Z} (-1)^{k+l} \bigl( \delta(k) - \delta(l) \bigr) q^{(k^2 + l^2)/3} x^l &= -x^{-1} q^{1/3} \frac{(q;q)_\infty}{(q^2;q^2)_\infty} \Theta(x;q^2) \Theta(x;q),\\
\sum_{k,l\in\Z} \bigl( \delta(k-1) - \delta(l) \bigr) q^{(k(k+1)/2 + 2l^2)/3} x^l &= -2 \frac{(q^2;q^2)_\infty}{(q;q)_\infty} \Theta(xq;q^2) \Theta(-xq^2;q^4).
\end{split}
\end{equation*}
\end{lemma}
\begin{proof}
Define for fixed $q$
\begin{equation*}
\begin{split}
f_L(x) &:= \sum_{k,l\in\Z} (-1)^{k+l} \bigl( \delta(k) - \delta(l) \bigr) q^{(k^2 + l^2)/3} x^l,\\
f_R(x) &:= x^{-1} \Theta(x;q^2) \Theta(x;q).
\end{split}
\end{equation*}
If we replace $l$ by $l+3$ in the definition of $f_L$ and use that $\Theta(x;q)= -x \Theta(qx;q)$, we see that
\begin{equation}\label{tr}
\begin{split}
f_L(x)&= -q^3 x^3 f_L(q^2 x),\\
f_R(x)&= -q^3 x^3 f_R(q^2 x).
\end{split}
\end{equation}
From the Jacobi triple product identity we see that $\Theta$ only has zeros if $x$ is an integer power of $q$. From this we see that $f_R$ has a double zero at $q^n$ if $n$ is even and a single zero at $q^n$ if $n$ is odd. $f_L$ has the same zeros:
to see that $f_L(1)=0$ replace $(k,l)$ by $(l,k)$ and that $f'_L(1)=0$ replace $l$ by $-l$; to see that $f_L(q) =0$ replace $l$ by $-l-3$. From \eqref{tr} we then get that $f_L$ also has a double zero at $q^n$ if $n$ is even and a single zero at $q^n$ if $n$ is odd. Hence $f_L / f_R$ has no poles and satisfies
\begin{equation*}
\frac{f_L}{f_R} (q^2 x) = \frac{f_L}{f_R} (x),
\end{equation*}
from which we get that $f_L / f_R$ is constant (as a function of $x$). To find this constant we consider the coefficients of $x^0$ in the expansion of both $f_L$ and $f_R$. The coefficient of $x^0$ in $f_L$ is
\begin{equation*}
\begin{split}
\sum_k (-1)^k \bigl( \delta(k) -1 \bigr) q^{k^2/3} &= - \sum_{k \not\equiv 0 \smod{3}} (-1)^k q^{k^2/3}= -2\sum_{k \equiv 1 \smod{3}} (-1)^k q^{k^2/3}\\
&= 2 q^{1/3} \sum_m (-1)^m q^{3m^2 +2m}= 2 q^{1/3} (q^6;q^6)_\infty (q;q^6)_\infty (q^5;q^6)_\infty \\
&= 2 q^{1/3} \frac{(q^6;q^6)_\infty (q;q^2)_\infty}{(q^3;q^6)_\infty}= 2 q^{1/3} \frac{(q;q)_\infty (q^6;q^6)_\infty^2}{(q^2;q^2)_\infty (q^3;q^3)_\infty}.
\end{split}
\end{equation*}
The coefficient of $x^0$ in $f_R$ is the coefficient of $x^1$ in 
\begin{equation*}
\sum_{n,m} (-1)^{n+m} q^{n(n-1) + m(m-1)/2} x^{n+m},
\end{equation*}
which is
\begin{equation*}
-\sum_n q^{3n(n-1)/2} = -2 (q^3;q^3)_\infty (-q^3;q^3)_\infty^2 = -2\frac{(q^6;q^6)_\infty^2}{(q^3;q^3)_\infty}.
\end{equation*}
Comparing these two we find that $f_L / f_R$ equals
\begin{equation*}
-q^{1/3} \frac{(q;q)_\infty}{(q^2;q^2)_\infty},
\end{equation*}
which gives the first identity. The proof of the second identity is very similar. Define for fixed $q$
\begin{equation*}
\begin{split}
g_L(x) &:= \sum_{k,l\in\Z} \bigl( \delta(k-1) - \delta(l) \bigr) q^{(k(k+1)/2 + 2l^2)/3} x^l,\\
g_R(x) &:= \Theta(xq;q^2) \Theta(-xq^2;q^4).
\end{split}
\end{equation*}
We can easily verify that both $g_L$ and $g_R$ satisfy
\begin{equation}\label{tr2}
g(x)= q^6 x^3 g(q^4 x).
\end{equation}
The zeros of $g_R$ are simple zeros at $x=q^n$ for $n$ odd and at $x=-q^n$ for $n\equiv 2 \pmod{4}$. To see that $g_L$ also has zeros in these points it suffices to show that $g_L(q)=0$ and $g_L(-q^2)=0$, by \eqref{tr2} and $g_L(x^{-1})=g_L(x)$. That $g_L(-q^2)=0$ follows directly if we replace $l$ by $-l-3$. Further 
\begin{equation*}
\begin{split}
g_L(q)&= \sum_{k,l\in\Z} \bigl( \delta(k-1) - \delta(l) \bigr) q^{(k(k+1)/2 + 2l^2+3l)/3}\\
&= q^{-1} \sum_{k\in\Z,\ m\ \text{odd}} \bigl( \delta(k-1) - \delta(m-1) \bigr) q^{(k(k+1)/2 + m(m+1)/2)/3}\\
&= \frac{1}{2} q^{-1} \sum_{k,m} \bigl( \delta(k-1) - \delta(m-1) \bigr) q^{(k(k+1)/2 + m(m+1)/2)/3}.
\end{split}
\end{equation*}
If we now replace $(k,m)$ by $(m,k)$ we get that $g_L(q)=0$. Since $g_L$ has zeros where $g_R$ has and from \eqref{tr2}, we get that $g_L / g_R$ is constant. Again we consider the coefficients of $x^0$ in both $g_L$ and $g_R$. For the coefficient in $g_L$ we have
\begin{equation*}
\begin{split}
\sum_k \bigl(\delta(k-1) -1 \bigr) q^{k(k+1)/6} &= -\sum_{k \not\equiv 1 \smod{3}} q^{k(k+1)/6} = -2 \sum_{k\equiv 0 \smod{3}} q^{k(k+1)/6}\\
&= -2 \sum_m q^{3m^2/2 +m/2} = -2 (q^3;q^3)_\infty (-q;q^3)_\infty (-q^2;q^3)_\infty \\
&= -2 \frac{(q^3;q^3)_\infty (-q;q)_\infty}{(-q^3;q^3)_\infty} = -2 \frac{(q^2;q^2)_\infty (q^3;q^3)_\infty^2}{(q;q)_\infty (q^6;q^6)_\infty}.
\end{split}
\end{equation*}
The coefficient of $x^0$ in $g_R$ is the coefficient of $x^0$ in
\begin{equation*}
\sum_{n,m} (-1)^n q^{n^2+2m^2}x^{n+m},
\end{equation*}
which is
\begin{equation*}
\sum_n (-1)^n q^{3n^2} = (q^6;q^6)_\infty (q^3;q^6)_\infty^2 = \frac{ (q^3;q^3)_\infty^2}{(q^6;q^6)_\infty}. 
\end{equation*}
Comparing these two we find that $g_L / g_R$ equals
\begin{equation*}
-2 \frac{(q^2;q^2)_\infty}{(q;q)_\infty},
\end{equation*}
which gives the second identity.
\end{proof}
We also need the following (slightly rewritten) result, which is Theorem 1.5 on page 646 in \cite{hick}.
\begin{lemma}\label{lem2}
Let $|q|<|x|<1$ and $|q|<|y|<1$ and $\Theta$ as \eqref{jac}. Then
\begin{equation*}
\sum_{r,s\in\Z} \rho_{r,s}\ x^r y^s q^{rs} = \frac{(q;q)_\infty^3 \Theta (xy;q)}{\Theta (x;q) \Theta (y;q)}.
\end{equation*}
\end{lemma}
We are now ready to prove Theorem \ref{theo}.
\begin{proof}
We let $C_{x^n y^m} [f]$ denote the coefficient of $x^n y^m$ in the expansion of $f$.

We start with the left hand side of the first equation
\begin{equation*}
\begin{split}
&\sum_{k,l,r,s\in\Z} \rho_{r,s} (-1)^{k+l+r+s} \bigl(\delta(k) - \delta (r) \bigr) \bigl(\delta(l) - \delta(s)\bigr)\ q^{(k^2+l^2+r^2+3rs+s^2 +3r+3s +1)/3}\\
&= C_{x^0 y^0} \left[ \sum_{k,l,m,n,r,s} \rho_{r,s} (-1)^{k+l+m+n} \bigl(\delta(k) - \delta (m) \bigr) \bigl(\delta(l) - \delta(n)\bigr)\ q^{(k^2+l^2+m^2+n^2)/3+rs+r+s+1/3} x^{m-r} y^{n-s}\right],
\end{split}
\end{equation*}
which by lemma \ref{lem1} and \ref{lem2} equals
\begin{equation*}
\begin{split}
C_{x^0 y^0} &\left[ -x^{-1} q^{1/3} \frac{(q;q)_\infty}{(q^2;q^2)_\infty} \Theta(x;q^2) \Theta(x;q) \cdot -y^{-1} q^{1/3} \frac{(q;q)_\infty}{(q^2;q^2)_\infty} \Theta(y;q^2) \Theta(y;q) \cdot q^{1/3} \frac{(q;q)_\infty^3 \Theta (x^{-1} y^{-1} q^2;q)}{\Theta (x^{-1} q;q) \Theta (y^{-1} q;q)}\right]\\
&= - \frac{(q;q)_\infty^5}{(q^2;q^2)_\infty^2} C_{x^0 y^0} \left[\Theta(x;q^2) \Theta(y;q^2) \Theta(x^{-1} y^{-1} q;q) \right]\\
&= - \frac{(q;q)_\infty^5}{(q^2;q^2)_\infty^2} C_{x^0 y^0} \left[\sum_{k,l,m} (-1)^{k+l+m} q^{k(k-1) +l(l-1) +m(m+1)/2} x^{k-m} y^{l-m} \right]\\
&= - \frac{(q;q)_\infty^5}{(q^2;q^2)_\infty^2} \sum_m (-1)^m q^{5m^2/2 -3m/2},
\end{split}
\end{equation*}
from which the first equation follows, using \eqref{theta}. The proofs of the other three are similar. Namely,
\begin{equation*}
\begin{split}
&\sum_{k,l,r,s\in\Z} \rho_{r,s} (-1)^{k+l+r+s} \bigl(\delta(k) - \delta (r-1) \bigr) \bigl(\delta(l) - \delta(s-1)\bigr)\ q^{(k^2+l^2+r^2+3rs+s^2 +r+s)/3} \\
&= C_{x^{-1} y^{-1}} \left[ \sum_{k,l,m,n,r,s} \rho_{r,s} (-1)^{k+l+m+n} \bigl(\delta(k) - \delta (m) \bigr) \bigl(\delta(l) - \delta(n)\bigr)\ q^{(k^2+l^2+m^2+n^2)/3+rs+r+s-2/3} x^{m-r} y^{n-s}\right]\\
&= - q^{-1} \frac{(q;q)_\infty^5}{(q^2;q^2)_\infty^2} C_{x^{-1} y^{-1}} \left[\sum_{k,l,m} (-1)^{k+l+m} q^{k(k-1) +l(l-1) +m(m+1)/2} x^{k-m} y^{l-m} \right]\\
&= - q^{-1} \frac{(q;q)_\infty^5}{(q^2;q^2)_\infty^2} \sum_{m} (-1)^m q^{5m^2/2-11 m/2+4}\\
&=\frac{(q;q)_\infty^5}{(q^2;q^2)_\infty^2} \sum_{m} (-1)^m q^{5m^2/2+ m/2},
\end{split}
\end{equation*}
where we replaced $m$ by $-m+1$ in the last step. Further
\begin{equation*}
\begin{split}
&\sum_{k,l,r,s\in\Z} \rho_{r,s} \bigl(\delta(k-1) - \delta (r) \bigr) \bigl(\delta(l-1) - \delta(s)\bigr)\  q^{(k(k+1)/2 + l(l+1)/2+2r^2+6rs +2s^2 +3r+3s)/3}\\
&= C_{x^0 y^0} \left[ \sum_{k,l,m,n,r,s} \rho_{r,s} \bigl(\delta(k-1) - \delta (m) \bigr) \bigl(\delta(l-1) - \delta(n)\bigr)\ q^{(k(k+1)/2 + l(l+1)/2+2m^2+2n^2)/3+2rs +r+s} x^{m-r} y^{n-s}\right],
\end{split}
\end{equation*}
which by lemma \ref{lem1} and \ref{lem2} equals
\begin{equation*}
\begin{split}
&= 4 \frac{(q^2;q^2)_\infty^5}{(q;q)_\infty^2} C_{x^0 y^0} \left[\Theta(xq;q^2) \Theta(-xq^2;q^4) \Theta(yq;q^2) \Theta(-yq^2;q^4) \cdot \frac{\Theta(x^{-1} y^{-1} q^2;q^2)}{\Theta(x^{-1} q;q^2)\Theta(y^{-1} q;q^2)}\right]\\
&= 4 \frac{(q^2;q^2)_\infty^5}{(q;q)_\infty^2} C_{x^0 y^0} \left[ \Theta(-xq^2;q^4)\Theta(-yq^2;q^4) \Theta(x^{-1} y^{-1} q^2;q^2)\right]\\
&= 4 \frac{(q^2;q^2)_\infty^5}{(q;q)_\infty^2} C_{x^0 y^0} \left[ \sum_{k,l,m} (-1)^m q^{2k^2 + 2l^2 + m^2+m} x^{k-m} y^{l-m} \right]= 4 \frac{(q^2;q^2)_\infty^5}{(q;q)_\infty^2} \sum_m (-1)^m q^{5m^2 +m},
\end{split}
\end{equation*}
from which the third equation follows, using \eqref{theta}. Finally we have
\begin{equation*}
\begin{split}
&\sum_{k,l,r,s\in\Z} \rho_{r,s} \bigl(\delta(k-1) - \delta (r+1) \bigr) \bigl(\delta(l-1) - \delta(s+1)\bigr)\  q^{(k(k+1)/2 + l(l+1)/2+2r^2+6rs +2s^2 +r+s-2)/3}\\
&= C_{x^0 y^0} \left[ \sum_{k,l,m,n,r,s} \rho_{r,s} \bigl(\delta(k-1) - \delta (m) \bigr) \bigl(\delta(l-1) - \delta(n)\bigr)\ q^{(k(k+1)/2 + l(l+1)/2+2m^2+2n^2)/3+2rs -r-s-2} x^{m-r-1} y^{n-s-1}\right]\\
&= 4 q^{-2} \frac{(q^2;q^2)_\infty^5}{(q;q)_\infty^2} C_{x^0 y^0} \left[x^{-1} y^{-1} \Theta(xq;q^2) \Theta(-xq^2;q^4) \Theta(yq;q^2) \Theta(-yq^2;q^4) \cdot \frac{\Theta(x^{-1} y^{-1} q^{-2};q^2)}{\Theta(x^{-1} q^{-1};q^2)\Theta(y^{-1} q^{-1};q^2)}\right]\\
&= 4 q^{-2} \frac{(q^2;q^2)_\infty^5}{(q;q)_\infty^2} C_{x^0 y^0} \left[ \Theta(-xq^2;q^4)\Theta(-yq^2;q^4) \Theta(x^{-1} y^{-1} q^{-2};q^2)\right]\\
&= 4 q^{-2} \frac{(q^2;q^2)_\infty^5}{(q;q)_\infty^2} C_{x^0 y^0} \left[ \sum_{k,l,m} (-1)^m q^{2k^2 + 2l^2 + m^2-3m} x^{k-m} y^{l-m} \right]= 4 \frac{(q^2;q^2)_\infty^5}{(q;q)_\infty^2} \sum_m (-1)^m q^{5m^2 -3m}.
\end{split}
\end{equation*}
\end{proof}


\begin{thebibliography}{99}
\bibitem{Choi1} Y.-S. Choi, {\it Tenth order mock theta functions in Ramanujan's Lost Notebook}, Invent.\ Math.\ \textbf{136} (1999), 497--569.
\bibitem{Choi2} Y.-S. Choi, {\it Tenth order mock theta functions in Ramanujan's Lost Notebook II}, Adv.\ Math.\ \textbf{156} (2000), 180--285.
\bibitem{Choi3} Y.-S. Choi, {\it Tenth order mock theta functions in Ramanujan's Lost Notebook III}, Proc.\ Lond.\ Math.\ Soc.\ (3) \textbf{94} (2007), 26–-52.
\bibitem{Choi4} Y.-S. Choi, {\it Tenth order mock theta functions in Ramanujan's Lost Notebook IV}, Trans.\ Amer.\ Math.\ Soc.\ \textbf{354} (2002), 705--733.
\bibitem{hick} D.R. Hickerson, {\it A proof of the mock theta conjectures}, Invent.\ Math.\ \textbf{94} (1988), 639--660.
\bibitem{Ono} K. Ono, {\it Unearthing the visions of a master: harmonic Maass forms and number theory}, Harvard-MIT Current Developments in Mathematics 2008, International Press, accepted for publication.
\bibitem{Ram} S. Ramanujan, {\it The lost notebook and other unpublished papers}, Narosa Publishing House, New Delhi, 1987.
\bibitem{Zag} D.B. Zagier, {\it Ramanujan's mock theta functions and their applications}, S\'eminaire Bourbaki, 2007-2008, no. 986.
\bibitem{zwe1} S.P. Zwegers, {\it Mock theta functions}, Ph.D. Thesis, Universiteit Utrecht, 2002.
\bibitem{zwe2} S.P. Zwegers, {\it Mock theta functions II: Indefinite theta functions of type $(r,1)$}, in preparation.
\end{thebibliography}
\end{document}